\documentclass[a4paper, 10pt, conference]{ieeeconf}      

\IEEEoverridecommandlockouts                              
\usepackage{amsmath}
\usepackage{amssymb}
\usepackage{epsfig}
\usepackage{psfrag}

\newcommand{\be}{\begin{equation}}
\newcommand{\ee}{\end{equation}}
\newcommand{\beano}{\begin{eqnarray*}}
\newcommand{\eeano}{\end{eqnarray*}}
\newcommand{\ba}{\begin{array}}
\newcommand{\ea}{\end{array}}

\newcommand{\Real}{\mathbb R}
\newcommand{\ra}{\rightarrow}
\newcommand{\norm}[1]{\|#1\|}
\newcommand{\bignorm}[1]{\left \| #1 \right \|}
\newcommand{\tr}{^{\rm T}}
\newcommand{\rform}[1]{(\ref{#1})}
\newcommand{\parmatrix}[1]{\begin{pmatrix} #1 \end{pmatrix}}

\newtheorem{theorem}[equation]{Theorem}
\newtheorem{remark}[equation]{Remark}
\newtheorem{assumption}{A\hspace{-1.5mm}}
\newtheorem{property}{P\hspace{-1.5mm}}

\title{\LARGE \bf
A Luenberger-style Observer for Robot Manipulators\\ with Position
Measurements}

\author{Fabio Celani
\thanks{This work was funded by the Research Council of Norway under the Centre of Excellence scheme.}
\thanks{The author is with the Centre for Ships and Ocean Structures (CeSOS),
NTNU Norwegian University of Science and Technology, Otto Nielsens
veg 10, NO-7491 Trondheim, Norway. {\tt\small celani@ntnu.no }.}%
}

\begin{document}

\maketitle
\thispagestyle{empty}
\pagestyle{empty}

\begin{abstract}

This paper presents a novel Luenberger-style observer for robot
manipulators with position measurements. Under the assumption that
the state evolutions that are to be observed have bounded
velocities, it is shown that the origin of the observation error
dynamics is globally exponentially stable and that the corresponding
convergence rate can be made arbitrarily high by increasing a gain
of the observer.

Comparisons and relations between the proposed observer and existing
observers are discussed.

The effectiveness of the result here presented is illustrated by a
simulation of the observer for the Pendubot, an underactuated
two-joint manipulator.

\end{abstract}

\section{INTRODUCTION}

Observers for nonlinear systems have been extensively studied in the
last decades (for example, see \cite{GauthierKupka01} and
\cite{NijmeijerFossen99}). In the specific case of robotic
manipulators it has been of great interest to design observers that
estimate velocities from position measurements. In fact, many
commercially-available robotic manipulators are not equipped with
velocity sensors; as a result, observers that estimate velocities
from position measurements have been instrumental in designing
feedback controllers that use only position measurements. The
literature on this topic is vast; refer to \cite{Xianetal04} for a
literature review. In some articles it is proved explicitly that, in
certain conditions, the state of the observer that is used to do
feedback from position measurements converges to the state of the
robotic manipulator even when the observer is in open-loop. Examples
of such articles are \cite{NicosiaTornambeValigi90},
\cite{NicosiaTomei90}, \cite{BerghuisNijmeijer93},
\cite{Khelfi_etal96}, and \cite{Xianetal04}. In addition,
\cite{CanudasSlotine91} and \cite{Khelfi_etal95} present observers
for robotic manipulators without employing them in a position
feedback loop. Moreover, rigid robotic systems are a subclass of
Euler-Lagrange systems, and intrinsic observers for Euler-Lagrange
systems with position measurements are proposed in
\cite{AghannanRouchon03} and \cite{AnisiHamberg05}.

This paper introduces a novel asymptotic observer for rigid robotic
manipulators with position measurements. The proposed observer is
Luenberger-style and is very simple in structure. Under the standard
and realistic assumption that the state evolutions that are to be
observed have bounded velocities with bounds known a priori, it is
shown that the origin of the observation error dynamics is globally
exponentially stable; in addition, it is demonstrated that the
corresponding convergence rate can be made arbitrarily high by
increasing a gain of the observer.

The observer here presented is compared with three similar existing
observers for robotic manipulators, and it is shown that in several
aspects the comparison is favorable to the proposed observer.

The rest of the paper is organized as follows; in section
\ref{main_res} the observer is presented; section
\ref{exist_results} discusses comparisons and relations between the
proposed observer and some existing results; in section \ref{sim_bb}
a simulation of the observer for the Pendubot is shown.

In what follows $\norm{\cdot}$ denotes the Eucledian norm of a
vector or matrix; moreover, given $v \in \Real^n$, $v_i$ denotes its
$i$-th component; $\Real_{>0}^n$ denotes the open interval $(0\
\infty)^n$.

\section{MAIN RESULT}\label{main_res}
The dynamic equation of a $n$-joint rigid robot can be expressed as
\be \label{robot} M(q) \ddot q + C(q,\dot q) \dot q + F \dot q +
g(q) = Hu \ee (see \cite{SciaviccoSiciliano00}). In \rform{robot},
$q \in \Real^n$ is the vector of joint positions, $u \in \Real^m$ is
the vector of control inputs, $M(q) \in \Real^{n \times n}$ is the
inertia matrix, $C(q,\dot q)$ is the Coriolis and centrifugal
matrix, $F$ is the diagonal matrix of viscous friction coefficients,
$g(q)$ is the gravity vector, and $H$ is the input matrix that
differs from the identity if the system is not fully actuated. $M$,
$C$, and $g$ are assumed to be continuously differentiable.

Denote by $Q \subseteq \Real^n$ the set of all feasible values of
$q$. $Q$ is typically determined by the mechanic structure of the
robot and/or by the way the robot operates. Note that if $q_i$ is
the joint position of a \emph{prismatic} joint, then $q_i$ ranges on
a bounded set because for physical reasons the relative linear
displacement of two links connected by a prismatic joint cannot grow
indefinitely. However, if $q_i$ is the joint position of a revolute
joint, then $q_i$ could range on an unbounded set since it can occur
that the links connected by a revolute joint can rotate indefinitely
with respect to each other. Thus, in the general case $Q$ is
unbounded.

The following properties of \rform{robot} are inherent to robot
dynamics (see \cite[p. 139]{CanudasdeWitFixotAstrom92}) and they
will be useful in the sequel  \vspace{2 mm}
\begin{property}
$\det (M(q)) \neq 0\ \ \forall q \in Q$
\end{property}
\vspace{2 mm}
\begin{property} \label{M_inv}
$\norm{M^{-1}(q)} \leq \ M_0\ \ \forall q \in Q$
\end{property}
\vspace{2 mm}
\begin{property} \label{coriolis}
$\forall i \in \{1,\ldots,n\}$, the $i$-th element of the vector
$C(q,\dot q) \dot q$ is equal to $\dot q\tr N_i(q) \dot q$ with
$N_i$ symmetric, continuously differentiable, and such that $\exists
\hat{N_i}>0$ that satisfies
\[ \norm{N_i(q)} \leq \hat{N_i}\ \ \forall q \in Q\;.\]
\end{property}
\vspace{2 mm}

It is assumed that the vector of joint positions \emph{$q$ is
measured}, but the vector of joint velocities \emph{$\dot q$ is not
measured}; then, \rform{robot} has the following state space
representation
\be \label{state_sp} \ba{rcl} \dot q &=& v \\
\dot v &=& -M^{-1}(q)(C(q,v) v + F v + g(q) -Hu) \\
 y &=& q \ea
\ee

This paper presents an observer for systems of the form
\rform{state_sp}; clearly, such observer is useful for estimating
the joint velocities $v$.

The state evolution of \rform{state_sp} that we want to observe
$(q(t),v(t))$ is assumed to be defined $\forall t \geq 0$ and with
bounded joint velocities, that is there exist $V_i \geq 0$
$i=1,\ldots,n$ such that \be \label{bounds} |v_i(t)| \leq V_i\ \
\forall t \geq 0\ \ \forall i \in \{1,\ldots,n\}\;,\ee and it is
assumed that the $V_i$'s are known a priori.

This assumption is definitively realistic. In fact, it is reasonable
to expect that the joint velocities of a robot will not exceed
certain a priori bounds that come from the mechanic limitations of
the robot and/or from the way the robot operates. Moreover, this
assumption is recurrent in the literature on observers for robotic
manipulators; for example, an equivalent assumption is made for
proving the convergence of the observers presented in
\cite{Xianetal04}, \cite{NicosiaTomei90},
\cite{BerghuisNijmeijer93}, \cite{Khelfi_etal96},
\cite{Khelfi_etal95} and \cite{AnisiHamberg05}.

Denote by $\hat q$ and $\hat v$ the estimates of $q$ and $v$
respectively; then, the proposed Luenberger-style observer is
defined by the following equations \be \label{observer} \ba{rcl}
\dot {\hat q} &=& \hat v - \theta \alpha (\hat q
-q)\\[1mm]
\dot {\hat v} &=& -M^{-1}(q)(C(q,\sigma_{V}(\hat v)) \sigma_{V}(\hat
v) + F \hat v + g(q)\\
&&-Hu) - \theta^2 \beta (\hat q-q) \ea \ee In \rform{observer}
$\alpha$, $\beta$, and $\theta$ are positive scalar design
parameters, $V=(V_1, \ldots, V_n)$ is the vector of the velocities
bounds, and $\sigma_V$ is a component-wise saturation function with
vector saturation level $V$; specifically, given $Y \in \Real^n$
such that $Y_i \geq~0$ $i=1,\ldots,n$, $\sigma_Y: \Real^n \ra
\Real^n$ is defined as follows \be \label{scal_sat}
(\sigma_Y(x))_i=\left\{ \ba{cl}
x_i & \mbox{if } |x_i| \leq Y_i\\
Y_i & \mbox{if } x_i > Y_i\\
-Y_i & \mbox{if } x_i < -Y_i \ea \right. \ee $i=1,\ldots,n$.

Observer \rform{observer} is obtained as follows. Make a copy of the
system \rform{state_sp} to be observed; add innovation terms to that
copy; saturate $\hat v$ in the Coriolis terms of the $\dot {\hat v}$
equation.

Note that the saturation level on each component of $\hat v$ is
chosen so that if the initial states of system \rform{state_sp} and
observer \rform{observer} are identical, then observer
\rform{observer} tracks exactly system \rform{state_sp}. In fact, if
$(q(t),v(t))$ is a state evolution of \rform{state_sp} corresponding
to a certain input $u(t)$, and it satisfies the bounds $|v_i(t)|
\leq V_i\ \ \forall t \geq 0\ \ \forall i \in \{1,\ldots,n\}$, then
$(q(t),v(t))$ is also the state evolution of the observer
\rform{observer} corresponding to the same input and to the initial
state $(\hat q(0), \hat v(0))=(q(0),v(0))$.
\\

The insertion of the saturation $\sigma_V$ in equations
\rform{observer} was inspired by \cite{ShimSonSeo01}. However, the
observer presented in \cite{ShimSonSeo01} applies to a class of
systems that does \emph{not} include systems of the type
\rform{state_sp}.\\

The following theorem states that observer \rform{observer} is
globally exponentially convergent with convergence rate arbitrarily fast.\\

\begin{theorem} \label{main_thm}
Let $(q(t),v(t))$ be the state evolution of \rform{state_sp}
corresponding to the input $u(t)$. Assume that $(q(t),v(t))$ is
defined $\forall t \geq 0$ and there exist $V_i \geq 0$
$i=1,\ldots,n$ such that \be \label{bounds2} |v_i(t)| \leq V_i \ \
\forall t \geq 0\ \ \forall i \in \{1,\ldots,n\}\;. \ee Then,
$\forall\ (\alpha, \beta, \gamma) \in \Real_{>0}^3\ \exists\
\theta^*>0$ such that if $\theta \geq \theta^*$
 the following property holds. $\exists k>0$ such that the state evolution $(\hat q(t),\hat v(t))$ of
\rform{observer} corresponding to the same input $u(t)$ and to
\textit{any} initial state $(\hat q(0),\hat v(0)) \in \Real^n \times
\Real^n$ is defined $\forall t \geq 0$ and satisfies \be
\label{upper_b} \left\| \parmatrix{\hat q(t)-q(t)\\ \hat
v(t)-v(t)} \right\| \leq k \left\| \parmatrix{\hat q(0)-q(0)\\
\hat v(0)-v(0)}\right\| e^{- \gamma t}\ \ \forall t \geq 0\;. \ee
\end{theorem}

\begin{proof}
To simplify the notation, let
\[
A(q,v) = C(q,v) v\;.
\]

Since \rform{bounds2} holds, in the rest of the proof regard
\be \label{gen_class_mod} \ba{rcl} \dot q &=& v \\
\dot v &=& -M^{-1}(q)(A(q,\sigma_V(v)) + F v + g(q) -Hu)\\ y &=& q
\ea \ee as the given system instead of \rform{state_sp}.

Fix $\alpha >0 \mbox{ and }  \beta >0$, and assume that $\theta >0$.
Similarly to \cite[p. 100]{GauthierKupka01} set
\[
\xi(t) = \frac{1}{\theta}(\hat q(t) -q(t))\ \ \ \zeta(t) =
\frac{1}{\theta^2}(\hat v(t) -v(t))\;.
\] Then, \be \label{error_dyn} \ba{rcl}
\parmatrix{\dot \xi(t) \\ \dot \zeta(t)} = \theta G \parmatrix{\xi(t) \\
\zeta(t)} + \parmatrix{0 \\ f(q(t),v(t),\zeta(t),\theta)} \ea \ee
where
\[
G = \parmatrix{- \alpha I & I \\ - \beta I & 0}
\] and
\begin{multline} \label{f} f(q,v,\zeta,\theta)=-M^{-1}(q)\\ \cdot \left\{ \left[
\frac{1}{\theta^2}(A(q,\sigma_V(v+\theta^2 \zeta)) -
A(q,\sigma_V(v))) \right] +F \zeta \right\} \;.
\end{multline} Note
that $G$ is Hurwitz since $\alpha$ and $\beta$ are positive. Let $S$
be the solution of the Lyapunov equation $G\tr S + S G = -I$, and
consider the candidate Lyapunov function for system
\rform{error_dyn}
\[
V(\xi,\zeta)= \parmatrix{\xi \\ \zeta}\tr S
\parmatrix{\xi \\ \zeta}\;.
\] Then
\begin{multline} \label{V_dot} \dot V (\xi,\zeta) \leq - \theta \left \| \parmatrix{\xi \\
\zeta} \right \|^2\\ + 2 \norm{S} \left \| \parmatrix{\xi \\
\zeta} \right \|
 \norm{f(q(t),v(t),\zeta,\theta)}\;.
\end{multline}

In order to find a proper upper bound for
$\norm{f(q(t),v(t),\zeta,\theta)}$, proceed as follows.

Using P\ref{coriolis} it follows that
\[
\frac{\partial A}{\partial v}(q,v)= 2 \left( \ba{c} v\tr N_1(q)\\
\vdots \\ v\tr N_n(q) \ea \right)\;.
\]  Let \be \label{barV} \bar V = \{v \in \Real^n
|\ |v_i| \leq V_i\ \ i=1,\ldots,n \}\;.\ee  Then, using
P\ref{coriolis} it follows that $\exists B >0$ such that \be
\label{L} \bignorm{ \frac{\partial A}{\partial v}(q,v)} \leq B\ \
\forall (q,v) \in Q \times \bar V \;. \ee Then, by \cite[Lemma
2]{ShimSonSeo01}
\begin{multline} \label{lipschitz} \norm{A(q,\sigma_V(v+\theta^2
\zeta)) - A(q,\sigma_V(v))} \leq \theta^2 B \norm{\zeta}\\ \forall
(q,v,\zeta,\theta) \in Q \times \Real^n \times \Real^n \times
\Real\;.
\end{multline}
Since $q(t) \in Q\ \ \forall t \geq 0$, letting $L=M_0
(B+\norm{F})$, from \rform{f}, \rform{lipschitz}, and P\ref{M_inv}
it follows that
\[
\norm{f(q(t),v(t),\zeta,\theta)}
 \leq L \norm{\zeta} \leq L \bignorm{\parmatrix{\xi \\ \zeta}}\;.
\]
Then, from \rform{V_dot}
\[
\dot V(\xi,\zeta) \leq -  \bignorm{\parmatrix{\xi \\
\zeta}}^2\left(\theta -  2\norm{S} L \right)\;.
\] As a result, if $\theta > 2 \norm{S} L$, the equilibrium at the origin of system \rform{error_dyn} is
globally exponentially stable. Standard calculations (see \cite[p.
154]{Khalil02}) show that the rate of the decaying exponential that
bounds from above $\norm{(\xi(t),\zeta(t))}$ is given by
\[ \frac{\theta - 2 \norm{S} L}{
2 \norm{S}}\;.\] As a result, to guarantee that this rate is greater
or equal than $\gamma$, it suffices to take \be \label{theta*}
\theta \geq \theta^* = 2 \norm{S} \left( \gamma + L \right) \ee
\end{proof}

\begin{remark} \label{roa}
Even though the state of the observer converges to the state of the
plant for any value of the initial state of the observer $(\hat
q(0),\hat v(0))$, it is enough to consider values of $(\hat
q(0),\hat v(0))$ with $\hat v(0) \in \bar V$ where $\bar V$ was
defined in \rform{barV}; in fact, it is known a priori that the
trajectory to be observed $(q(t),v(t))$ is such that $v(t) \in \bar
V\ \forall t \geq 0$. Moreover, since $q$ is measurable, it should
be feasible to set $\hat q(0) \approx q(0)$; consequently, in
practice it is enough to guarantee that the observer converges when
$(\hat q(0),\hat v(0)) \in \{(\hat q(0),\hat v(0))|\ \norm{\hat
q(0)-q(0)} < \epsilon, \hat v(0) \in \bar V \}$ where $\epsilon >0$
is a parameter whose value depends on the accuracy of the position
sensors.
\end{remark}

\section{Comparisons and Relations with Existing Results} \label{exist_results}

The proposed observer is derived under assumptions equivalent to
those for the observers in \cite{NicosiaTomei90},
\cite{BerghuisNijmeijer93}, and \cite{Khelfi_etal96} used in
open-loop; moreover, those observers and the one here proposed
present similar convergence properties. However, the observer here
introduced compares favorably to those in \cite{NicosiaTomei90},
\cite{BerghuisNijmeijer93}, and \cite{Khelfi_etal96} because it is
simpler in structure and consequently easier to implement. Indeed,
the proposed observer is a plain Luenberger-style observer with a
saturation on some of the $\hat v$ terms. Note that the observer in
\cite{Khelfi_etal96} has the advantage over the observer here
proposed  of being of reduced order; however, the structure of the
former is quite complicated and, as pointed out in
\cite{BonaIndri98}, the procedure to choose its parameters is quite
complex.

From a mathematical point of view, the observer here introduced
compares favorably with those in \cite{NicosiaTomei90} and
\cite{BerghuisNijmeijer93}. In fact, in the proposed observer the
error dynamics have the origin that is globally asymptotically
stable; on the other hand, the origin of the error dynamics is only
semiglobally stabilized in the case of the observers in
\cite{NicosiaTomei90} and \cite{BerghuisNijmeijer93}. However,
taking into account the considerations in Remark \ref{roa}, it
follows that achieving global rather than semiglobal convergence
might not be relevant from a practical point of view. An additional
point in favor of the observer here presented with respect to
observers in \cite{NicosiaTomei90}, \cite{BerghuisNijmeijer93}, and
\cite{Khelfi_etal96}, is that the proof of its convergence is
simpler.

The proposed observer is related to the nonlinear observer
introduced in \cite{TarguiFarzaHammouri02} as discussed in the rest
of the section.

The nonlinear observer presented in \cite{TarguiFarzaHammouri02}
applies to a certain class of multi-output nonlinear systems that
includes systems of the form \be \label{mons} \ba{rcl}
\dot q_1 &=& f^1_1(u,q,v_1)\\[1 mm]
\dot q_2 &=& f^1_2(u,q,v_1,v_2)\\
& \vdots &\\
\dot q_n &=& f^1_n(u,q,v)\\[1 mm]
\dot v_1 &=& f^2_1(u,q,v)\\
& \vdots &\\
\dot v_n &=& f^2_n(u,q,v)\\ \\
y &=& q
 \ea
\ee where $q_i,v_i \in \Real\ i=1,\ldots,n,\ u \in \Real^m$, and it
is assumed that
\begin{assumption} \label{u_compact}
$u(t) \in U$ a compact subset of $\Real^m$.
\end{assumption}
\begin{assumption} \label{gl_lipschitz}
$\forall (k,i) \in \{1,2\} \times \{1,\ldots,n\}\ f_i^k \in C^1$ and
$f_i^k$ is globally Lipschitz with respect to $(q,v)$ uniformly with
respect to $u \in U$.
\end{assumption}
\begin{assumption}\label{bdd}
$\exists\ \ 0< \bar \alpha < \bar \beta$ such that \begin{multline*}
0 < \bar \alpha \leq \frac{\partial f_i^1}{\partial v_i}(u,q,v) \leq
\bar \beta \ \ \forall (u,q,v) \in U \times \Real^n \times \Real^n\\
\forall i \in \{1,\ldots,n\}.
\end{multline*}
\end{assumption} Let $\alpha$, $\beta$, and $\theta$ be scalars;
then, in \cite{TarguiFarzaHammouri02} the following Luenberger-style
observer for system \rform{mons} is proposed \be \label{obs_mons}
\ba{rcl} \dot {\hat q}_1 &=& f^1_1(u,\hat q,\hat v_1)- \theta \alpha
(\hat q
-q)\\[1 mm]
\dot {\hat q}_2 &=& f^1_2(u,\hat q,\hat v_1,\hat v_2)- \theta \alpha
(\hat q
-q)\\
& \vdots &\\
\dot {\hat q}_n &=& f^1_n(u,\hat q,\hat v)- \theta \alpha (\hat q
-q)\\[1 mm]
\dot {\hat v}_1 &=& f^2_1(u,\hat q,\hat v)- \theta^2 \beta (\hat q-q)\\
& \vdots &\\
\dot {\hat v}_n &=& f^2_n(u,\hat q,\hat v)- \theta^2 \beta (\hat
q-q)\;.\\ \ea \ee

In \cite{TarguiFarzaHammouri02} it is proved that, $\forall\
(\alpha, \beta) \in \Real_{>0}^2\ \exists\ \theta^*>0$ such that if
$\theta > \theta^*$ observer \rform{obs_mons} is globally
exponentially convergent.

An alternative observer for \rform{mons} is given by \be
\label{obs_mons_alt} \ba{rcl} \dot {\hat q}_1 &=& f^1_1(u,q,\hat
v_1)- \theta \alpha (\hat q
-q)\\[1 mm]
\dot {\hat q}_2 &=& f^1_2(u,q,\hat v_1,\hat v_2)- \theta \alpha
(\hat q
-q)\\
& \vdots &\\
\dot {\hat q}_n &=& f^1_n(u,q,\hat v)- \theta \alpha (\hat q
-q)\\[1 mm]
\dot {\hat v}_1 &=& f^2_1(u,q,\hat v)- \theta^2 \beta (\hat q-q)\\
& \vdots &\\
\dot {\hat v}_n &=& f^2_n(u,q,\hat v)- \theta^2 \beta (\hat
q-q)\;.\\ \ea \ee It can be easily proved again that $\forall\
(\alpha, \beta) \in \Real_{>0}^2\ \exists\ \theta^*>0$ such that if
$\theta
> \theta^*$ observer \rform{obs_mons_alt} is globally exponentially
convergent. The advantage of \rform{obs_mons_alt} over
\rform{obs_mons} is that in order to prove convergence, assumption
A\ref{gl_lipschitz} can be replaced by the following weaker
assumption\\

\begin{assumption} \label{gl_lipschitz_alt}
$\forall (k,i) \in \{1,2\} \times \{1,\ldots,n\}\ f_i^k \in C^1$
and, denoting with $Q \subseteq \Real^n$ the set $q$ ranges on,
$f_i^k$ is globally Lipschitz with respect to $v$ uniformly with
respect to $(q,u) \in Q \times U$.\\
\end{assumption}

Clearly, system \rform{state_sp} is of the type \rform{mons} and it
satisfies assumption A\ref{bdd}; however, in general, it does not
satisfy assumption A\ref{gl_lipschitz_alt}. On the other hand, it is
assumed for system \rform{state_sp} that the velocity $v(t)$ of the
state evolutions to be observed satisfy \rform{bounds2}.
Consequently, as said before, system
\be \label{state_sp_satv} \ba{rcl} \dot q &=& v \\
\dot v &=& -M^{-1}(q)(C(q,\sigma_V(v)) \sigma_V(v) + F v + g(q)\\
&& -Hu) \\
 y &=& q \ea
\ee  can be regarded as the given system instead of
\rform{state_sp}. Note that \emph{the proposed observer
\rform{observer} coincides with observer \rform{obs_mons_alt}
instanced for system \rform{state_sp_satv}}. Convergence of
\rform{observer} can be justified as follows. From what stated
before, it follows that if A\ref{u_compact}, A\ref{bdd}, and
A\ref{gl_lipschitz_alt} hold for \rform{state_sp_satv} then
convergence is achieved. Note that \rform{state_sp_satv} satisfies
A\ref{bdd} and, as shown in the proof of Theorem \ref{main_thm},
using properties P\ref{M_inv}, P\ref{coriolis}, and \cite[Lemma
2]{ShimSonSeo01}, it follows that \rform{state_sp_satv} satisfies
A\ref{gl_lipschitz_alt}, too. Moreover, note that assumption
A\ref{u_compact} is not needed to prove convergence because in
equations \rform{state_sp_satv} $u$ enters only through the additive
term $M^{-1}(q)Hu$ which does not depend on the unmeasured variable
$v$; as a result, the dynamics of the observation error are
independent of $u$.

\section{Simulation of the Observer for the Pendubot}\label{sim_bb} 
\begin{figure}
\begin{center}
  \centerline{\epsfig{file=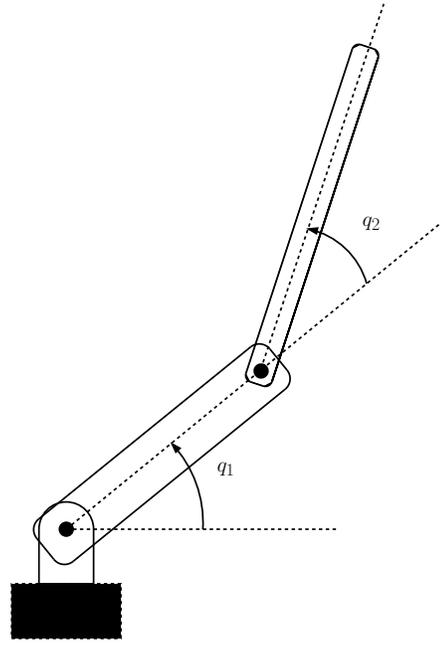, height=8.5 cm} }
  \caption{\label{pend_sketch} Sketch of the Pendubot.}
\end{center}
\end{figure}

The
effectiveness of the proposed Luenberger-style observer is here
illustrated by a simulation of the observer for the Pendubot, an
underactuated two-joint manipulator moving in a vertical plane (see
\cite{SpongBlock95}). A sketch of the Pendubot is shown in Fig.
\ref{pend_sketch}.
The first joint (shoulder) is actuated, whereas the second joint
(elbow) is not. Both joints are equipped with position sensors
(encoders), but no velocity sensors are present. Consequently, it is
of interest to design an observer for the Pendubot that estimates
the joint-velocities so that the observer could be included in a
feedback controller that uses only position measurements.

Define the joint variables $q_1$ and $q_2$ as shown Fig.
\ref{pend_sketch}, and let $u$ be the voltage input of the actuator.
A dynamic model of the Pendubot can be found in \cite{ZhangTarn02}
and is given by \rform{robot} with $q=(q_1\ q_2)\tr$,
\[
M(q) = \left( \ba{cc} \pi_1 + \pi_2 + 2 \pi_3 \cos(q_2) & \pi_2+\pi_3\cos(q_2)\\
\pi_2+\pi_3\cos(q_2) & \pi_2 \ea\right) \]
\begin{multline*}
C(q,\dot q) \\
= \left(\ba{cc} -\pi_3\sin(q_2)\dot q_2 &
-\pi_3\sin(q_2)\dot q_2 -\pi_3\sin(q_2)\dot q_1 \\
\pi_3\sin(q_2)\dot q_1  & 0\ea \right)
\end{multline*}
\[g(q)=\left(\ba{c} \pi_4 g_0 \cos(q_1)+\pi_5g_o\cos(q_1+q_2)\\ \pi_5g_0\cos(q_1+q_2)\ea\right)\]
\[
F_v=\left(\ba{cc} 0 & 0\\ 0 & 0 \ea \right)\ \ \ H = \left(\ba{cc}
1& 0\\ 0 & 0 \ea \right)\;
\] where
\beano
\pi_1 &=& 0.0308\text{ Vs}^2/\text{rad}\\
\pi_2 &=& 0.0106\text{ Vs}^2/\text{rad}\\
\pi_3 &=& 0.0095\text{ Vs}^2/\text{rad}\\
\pi_4 &=& 0.2086\text{ Vs}^2 \text{/m}\\
\pi_5 &=& 0.0630\text{ Vs}^2 \text{/m}\\
g_0 &=& 9.81 \text{ m/s}^2\;.\eeano

Assume that the Pendubot operates so that the angular velocities
$v=\dot q$
 do not exceed the following bounds
\[
|v_1(t)| \leq 10\text{ rad/s}\ \ |v_2(t)| \leq 10\text{ rad/s}\
\forall t \geq 0\;.
\]
The design parameters $\alpha$ and $\beta$ of observer
\rform{observer} are set as $\alpha = \beta = 1$; then, referring to
magnitudes introduced in the proof of Theorem \ref{main_thm}, it
follows that $\norm{S} = 1.81$ and that $L$ can be set equal to
54.01; consequently, the minimum value of the gain $\theta$ that
guarantees global exponential stability of the origin of the error
dynamics is $\theta^*=195$. Set $\theta = 200$ so that the norm of
the observation error will converge to $0$ globally, and it will be
bounded by an exponential as in \rform{upper_b} with $\gamma=1.27$.

Choose the following initial state for the Pendubot
\[ \left( \ba{c} q_1^0\\ q_2^0\\ v_1^0\\ v_2^0 \ea \right) =
\left( \ba{c} -\pi/2\\ 0 \\ 0\\
0\ea \right)\;,
\] which corresponds to the lower equilibrium, and apply the control $u = 1.5 \sin(100 t)$ that maintains the Pendubot
in oscillation about the lower equilibrium.

\begin{figure}
\begin{center}
\centerline{\epsfig{file=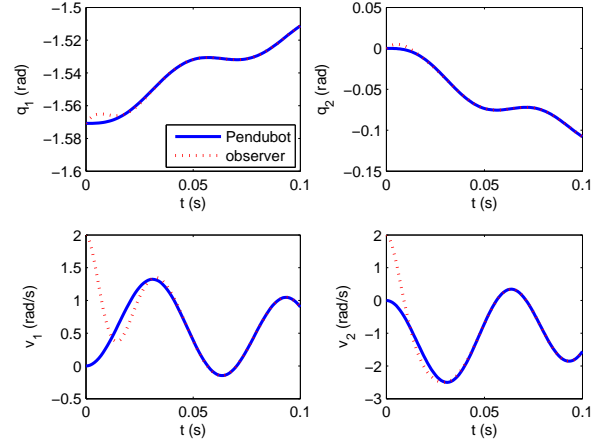, width=8.5 cm} }
\caption{\label{plot_pend_obs} Simulation of the Pendubot and its
observer.}
\end{center}
\end{figure}

Taking into account that $q$ is measured, choose for the observer
the initial state
\[ \left( \ba{c} \hat q_1^0\\ \hat q_2^0\\ \hat v_1^0\\ \hat v_2^0 \ea \right) =
\left( \ba{c} -\pi/2\\ 0\\ 2\\
2 \ea \right)\;.
\] The corresponding state evolutions of the Pendubot and of the
observer, plotted in Fig.~\ref{plot_pend_obs}, show that the
observer is convergent.

\section{Conclusions}
In this paper a novel asymptotic Luenberger-style observer for robot
systems with position measurements has been presented. The observer
is very simple in structure; it has been proven that the dynamics of
the observation error have a globally exponentially stable origin
with convergence rate that can be made arbitrarily high by
increasing a gain of the observer.

The proposed observer compares favorably with some existing
observers for robot manipulators derived under equivalent
assumptions; its relation to a previous nonlinear observer has been
discussed.

A simulation of the proposed observer for the Pendubot has been
included to illustrate its effectiveness.

It will be topic of future research to investigate if this observer
can lead to interesting results in the area of control of robot
manipulators via position feedback .

\section{Acknowledgment}
The author is grateful to Claudio De Persis and Dennis Lucarelli for
their valuable suggestions.


\bibliographystyle{IEEEtran}
\bibliography{IEEEabrv,arxiv_rob_obs_05}

\begin{thebibliography}{10}
\providecommand{\url}[1]{#1}
\csname url@rmstyle\endcsname
\providecommand{\newblock}{\relax}
\providecommand{\bibinfo}[2]{#2}
\providecommand\BIBentrySTDinterwordspacing{\spaceskip=0pt\relax}
\providecommand\BIBentryALTinterwordstretchfactor{4}
\providecommand\BIBentryALTinterwordspacing{\spaceskip=\fontdimen2\font plus
\BIBentryALTinterwordstretchfactor\fontdimen3\font minus
  \fontdimen4\font\relax}
\providecommand\BIBforeignlanguage[2]{{%
\expandafter\ifx\csname l@#1\endcsname\relax
\typeout{** WARNING: IEEEtran.bst: No hyphenation pattern has been}%
\typeout{** loaded for the language `#1'. Using the pattern for}%
\typeout{** the default language instead.}%
\else
\language=\csname l@#1\endcsname
\fi
#2}}

\bibitem{GauthierKupka01}
J.~Gauthier and I.~Kupka, \emph{Deterministic Observation Theory and
  Applications}.\hskip 1em plus 0.5em minus 0.4em\relax Cambridge: Cambridge
  University Press, 2001.

\bibitem{NijmeijerFossen99}
H.~Nijmeijer and T.~Fossen, Eds., \emph{New Directions in Nonlinear Observer
  Design}.\hskip 1em plus 0.5em minus 0.4em\relax London: Springer-Verlag,
  1999.

\bibitem{Xianetal04}
B.~Xian, M.~de~Queiroz, D.~Dawson, and M.~McIntyre, ``A discontinuous output
  feedback controller and velocity observer for nonlinear mechanical systems,''
  \emph{Automatica}, vol.~40, no.~4, pp. 695--700, 2004.

\bibitem{NicosiaTornambeValigi90}
S.~Nicosia, A.~Tornambe, and P.~Valigi, ``Experimental results in state
  estimation of industrial robots,'' in \emph{Proceedings of the 29th IEEE
  Conference on Decision and Control}, December 1990, pp. 360--5.

\bibitem{NicosiaTomei90}
S.~Nicosia and P.~Tomei, ``Robot control by using only joint position
  measurements,'' \emph{{IEEE} Trans. Automat. Contr.}, vol.~35, no.~9, pp.
  1058--61, 1990.

\bibitem{BerghuisNijmeijer93}
H.~Berghuis and H.~Nijmeijer, ``A passivity approach to controller-observer
  design for robots,'' \emph{{IEEE} Trans. Robot. Automat.}, vol.~9, no.~6, pp.
  740--34, 1993.

\bibitem{Khelfi_etal96}
M.~Khelfi, M.~Zasadzinski, H.~Rafaralahy, E.~Richard, and M.~Darouach,
  ``Reduced-order observer-based point-to-point and trajectory controllers for
  robot manipulators,'' \emph{Control Engineering Practice}, vol.~4, no.~7, pp.
  991--1000, 1996.

\bibitem{CanudasSlotine91}
C.~Canudas~de Wit and J.-J. Slotine, ``Sliding observers for robot
  manipulators,'' \emph{Automatica}, vol.~27, no.~5, pp. 859--64, 1991.

\bibitem{Khelfi_etal95}
M.~Khelfi, M.~Zasadzinski, M.~Darouach, and E.~Richard, ``Reduced-order
  ${H}_{\infty}$ observer for robot manipulators,'' in \emph{Proc. {IEEE}
  Conference on Decision and Control}, New Orleans, LA, USA, Dec. 1995, pp.
  1011--1013.

\bibitem{AghannanRouchon03}
N.~Aghannan and P.~Rouchon, ``An intrinsic observer for a class of {L}agrangian
  systems,'' \emph{{IEEE} Trans. Automat. Contr.}, vol.~48, no.~6, pp. 936--45,
  2003.

\bibitem{AnisiHamberg05}
D.~Anisi and J.~Hamberg, ``Riemannian observers for {E}uler-{L}agrange
  systems,'' in \emph{Proc. 16th IFAC World Congress, to appear}.

\bibitem{SciaviccoSiciliano00}
L.~Sciavicco and B.~Siciliano, \emph{Modelling and Control of Robot
  Manipulators}.\hskip 1em plus 0.5em minus 0.4em\relax London: Springer, 2000.

\bibitem{CanudasdeWitFixotAstrom92}
C.~Canudas~de Wit, N.~Fixot, and K.~J. Astrom, ``Trajectory tracking in robot
  manipulators via nonlinear estimated state feedback,'' \emph{{IEEE} Trans.
  Robot. Automat.}, vol.~8, no.~1, pp. 138--144, 1992.

\bibitem{ShimSonSeo01}
H.~Shim, Y.~Son, and J.~Seo, ``Semi-global observer for multi-output nonlinear
  systems,'' \emph{Systems and Control Letters}, vol.~42, no.~3, pp. 233--44,
  2001.

\bibitem{Khalil02}
H.~Khalil, \emph{Nonlinear Systems}.\hskip 1em plus 0.5em minus 0.4em\relax
  Upper Saddle, NJ: Prentice Hall, 2002.

\bibitem{BonaIndri98}
B.~Bona and M.~Indri, ``Analysis and implementation of observers for robotic
  manipulators,'' in \emph{Proc. {IEEE} Conference on Robotics and Automation},
  Leuven, Belgium, May 1998, pp. 3006--11.

\bibitem{TarguiFarzaHammouri02}
B.~Targui, M.~Farza, and H.~Hammouri, ``Constant-gain observer for a class of
  multi-output nonlinear systems,'' \emph{Applied Mathematics Letters},
  vol.~15, no.~6, pp. 709--20, 2002.

\bibitem{SpongBlock95}
M.~Spong and D.~Block, ``The {P}endubot: a mechatronic system for control
  research and education,'' in \emph{Proc. {IEEE} Conference on Decision and
  Control}, New Orleans, LA, USA, Dec. 1995, pp. 555--6.

\bibitem{ZhangTarn02}
M.~Zhang and T.-J. Tarn, ``Hybrid control of the {P}endubot,''
  \emph{{IEEE/ASME} Trans. Mechatron.}, vol.~7, no.~1, pp. 79--86, 2002.

\end{thebibliography}

\end{document}